\documentclass[12pt]{article}
\usepackage[latin1]{inputenc}
\usepackage{a4wide,amssymb}
\usepackage[francais]{babel}
\usepackage[T1]{fontenc}
\usepackage{pifont}
\usepackage{amsmath,amscd}
\usepackage{graphics}
\usepackage{pb-diagram}
\usepackage{pb-diagram,pb-xy}
\usepackage{fancyheadings}
\usepackage{amsmath,amssymb,nopageno}
\usepackage[all]{xy}
\newtheorem{thm}{Th\'{e}or\`{e}me}
\newtheorem{deft}{D\'{e}finition}
\newtheorem{proposition}{Proposition}
\newtheorem{rmq}{Remarque}

\newtheorem{lemme}{Lemme}

\newtheorem{corr}{Corollaire}
\begin{document}

%\tableofcontents
%\newpage

%%%%%%%%%%%%%%%%%%%      Intro        %%%%%%%%%%%%%%%%%%%%%%%%%%%%%%%%%%%%%%%

%%%%%%%%%%%%%%%%%%%%%%%%%%%%%%%%%%%%%%%%%%%%%%%%%%%%%%%%%%%%%%%%%%%%%%%%%%%%%%%%%%%%%%%%%%%%%%%%%%%%%%%%%%%%%%%%%%%%%%%%%%%%%%%%%%%%%%%%%%%%%ù

\title{G\'eom\'etricit\'e artinienne de l'$\infty$-champs des \'el\'ements de Maurer-Cartan}

\author{Brahim Benzeghli\thanks{Ce papier a b\'en\'efici\'e d'une aide
de l'Agence Nationale de la Recherche portant
la r\'ef\'erence ANR-09-BLAN-0151-02 (HODAG)}}

%\begin{abstract}
 %L'abstract va ici
%\end{abstract}

\maketitle

{\small
\noindent
{\bf R\'esum\'e:}
Dans cet article, en se basant sur les m\^emes techniques que dans \cite{BENZ08} pour la construction de la carte formellement lisse $V \to Perf$  entre la vari\'et\'e des complexes
et le $\infty$-champs d'Artin $Perf$, on construira explicitement un nouveau $\infty$-champs d'Artin des \'el\'ements de Maurer-Cartan d'une dg-cat\'egorie $\mathcal{P}$, 
qu'on note par $\mathcal{MC_P}$, avec une carte $V \to \mathcal{MC_P}$ qu'on montrera formellement lisse.  
}
\\
\\
{\small
\noindent
{\bf Abstract:}
In this paper, based on the same techniques as in \cite{BENZ08} for the construction of a formally smooth map  $V \to Perf$ from the variety of complexes  
to  the Artin $\infty$-stack $Perf$,  we explicitly construct a new Artin $ \infty $-stack of Maurer-Cartan  elements
of a dg-category $ \mathcal {P} $ denoted by $\mathcal{MC_P} $, with a map $ V \to \mathcal{MC_P} $ which we show to be formally smooth.
}
\section*{Introduction}

Dans  \cite{BENZ08} on a construit une carte explicite $V \to Perf$ pour l'$ \infty$-champs des complexes parfaits, o\`u $V$ \'etait le sch\'ema de Buchsbaum-Eisenbud \cite{BUCH1}, \cite{BUCH2}, \cite{BRUN}, \cite{HUNE}, \cite{KEMP}, \cite{MASS},  \cite{TRIV} et \cite{YOSH} 
qui param\'etrise les diff\'erentiels $d$ avec $d^2=0$ sur une suite de fibr\'es vectoriels triviaux.
\\
Le but principal du \cite{BENZ08} \'etait de montrer la lisset\'e formelle du morphisme  $V \to Perf$, apr\'es avoir explicit\'e l'$ \infty-$champs d'Artin $Perf$. 
\\

On veut g\'en\'eraliser ce r\'esultat pour un autre champs. Pour cela on fixera  une $dg$-cat\'egorie $k$-lin\'eaire $\mathcal P$ qui satisfait aux hypoth\`eses suivants:
\begin{itemize}
 \item L'ensemble des objets $Ob(\mathcal P)$ est fini.
\item Pour tout $E,F \in Ob(\mathcal P)$, pour tout $i \in \mathbb Z, \quad \mathcal P^i(E,F) $ est un $k$-espace vectoriel de dimension fini.
\item Il existe un indice $n > 0$ tel que pour tout $i< -n$, le $k$-espace vectoriel  $\mathcal P^i(E,F)=0$. 
\end{itemize}

On peut d\'efinir une $(\infty,1)-$cat\'egorie $MC({\mathcal P})$ dont les objets sont les couples $(E, \eta)$ o\`u $E$ est un objet de $\mathcal{P}$ et $\eta$ est un \'el\'ement de 
Maurer-Cartan dans $\mathcal{P}^1(E,E)$.
\\

On d\'efinit l'$\infty$-champs $\mathcal{MC_P}$ comme le $\infty$-champs associ\'e \`a l'$\infty$-pr\'echamps $\mathcal{MC_P}^{es}$ qui \`a une $k$-alg\`ebre $B$ associe l'int\'erieur de $MC({\mathcal P} \otimes_kB)$.
Cet int\'erieur peut \^etre vu comme un ensemble simplicial de $Kan$ ou une quasi-cat\'egorie $\mathcal{MC_P}^{es}(B) $.
\\

Pour construire une carte recouvrant l'$\infty$-champs $\mathcal{MC_P} $ on construira un foncteur
 $$
 V_E : AlgCom_k \to \mathcal Ens
 $$
  qui associe \`a chaque $B \in AlgCom_k$ son image $V_E(B)$ l'ensemble des \'el\'ements de Maurer-Cartan dans $\mathcal{P}^1(E,E)\otimes_kB $. Ce foncteur est repr\'esentable par un sch\'ema affine.
  On proc\`ede en choisissant des $k$-bases dans $\mathcal{P}^1$ et $\mathcal{P}^2$ donn\'ees et en suivant les m\^emes d\'emarches que dans \cite{BENZ08};  $V_E$ est le pr\'eimage de l'origine d'une application entre espaces affines

\[ \begin{array}{cccc}
courb : &\mathcal{P}^1_{sch}(E,E)&  \to  & \mathcal{P}^2_{sch}(E,E)\\
   {}   & \eta & \mapsto & d(\eta) + \eta  \eta .
\end{array}  \]

La d\'etermination du $V_E$ nous permettra d'avoir une carte 
$$
V_E \to \mathcal{MC_P}.
$$
On d\'emontre qu'elle est formellement lisse. 
Pour \c{c}a on montre que pour un anneau commutatif $B$ et si $I$ un id\'eal de carr\'e nul dans $B$, alors: pour tout couple de morphismes
$$
(F, \zeta_I) \quad :  \quad \overset{  \alpha_I  }{\underset{  a_I }{ \rightleftarrows}} \quad : \quad  ( E,\eta_I := \eta \otimes _B B/I )
$$
tels que $ a_I$ et $ \alpha_I$ sont inverses dans $H^0(\mathcal{MC}(\mathcal{P} \otimes_k B/I)) $,
alors il existe un rel\`evement $\zeta$ pour $\zeta_I$ et $\alpha$ pour $\alpha_I$ tel que $d(\zeta) + \zeta^2 = 0$ et $d_{\zeta \eta}(\alpha) =0$.

L'existence d'une carte permet de d\'eduire que $\mathcal{MC_P}$ est un $(n+1)$-champs d'Artin.

 %PAGE 01 

\section*{  Des dg-Cat\'egories vers les $( \infty,1)$-Cat\'egories  }
\begin{deft}
 Soit $B$ un anneau commutatif. Une $dg$-cat\'{e}gorie $\mathcal{P}$ sur $B$ est donn\'{e}e par:
\begin{itemize}
 \item[i)] $ob( \mathcal{P})$ est un ensemble.
\item[ii)] $\forall E,F \in ob( \mathcal{P}) $, on a un complexe de $B$-modules $ \mathcal{P}^{ \cdot} (E,F)$ tel que
$$
( \mathcal{P}^{ \cdot} (E,F),d) = \{  (\mathcal{P}^{ i} (E,F), d^i), \quad d^i:  \mathcal{P}^{ i} (E,F) \rightarrow  \mathcal{P}^{ i+1} (E,F)  \}.
$$
autrement dit :
$$
(\mathcal{P}^{ \cdot} (E,F), d) = \quad ... \rightarrow \mathcal{P}^{i-1} (E,F) \xrightarrow{d^{i-1}}  \mathcal{P}^{ i} (E,F) \xrightarrow{d^i}  \mathcal{P}^{ i+1} (E,F) \rightarrow ...
$$

\item[iii)] $ 1 \in  \mathcal{P}^{ 0} (E,F)$  et  $ d(1)=0$.
\item[iv)] $ \forall \alpha \in  \mathcal{P}^{ i} (E,F)$, $ \forall \beta \in  \mathcal{P}^{ j} (F,G) $ alors $ \beta . \alpha \in  \mathcal{P}^{ i+j} (E,G)$ v\'{e}rifiant:
\begin{itemize}
\item[$\bullet$] L'application $ \alpha, \beta \mapsto \beta . \alpha $ est $B$-bilin\'eaire.
 \item[$\bullet$] $ \gamma ( \beta \alpha ) = (\gamma \beta  ) \alpha \quad \forall \alpha \in  \mathcal{P}^{ i} (E,F),  \beta \in  \mathcal{P}^{ j} (F,G), \gamma \in  \mathcal{P}^{k} (G,H)  $ .
\item[$\bullet$] $ d( \beta \alpha) = d(\beta) \alpha + (-1)^{ \rvert \beta \rvert } \beta d(\alpha)$.

\end{itemize}

\end{itemize}
\end{deft}

 Soient $ \mathcal{P}^{ \cdot} (X,Y)  $, $\mathcal{P}^{ \cdot} (Y,Z)$  et $\mathcal{P}^{ \cdot} (X,Z)$ trois complexes, on d\'efinit un morphisme
\begin{eqnarray*}
\mathcal{P}^{ \cdot} (Y,Z) \otimes    \mathcal{P}^{ \cdot} (X,Y)   & \rightarrow &  \mathcal{P}^{ \cdot} (X,Z) \\
\sum   \beta_i \otimes  \alpha_i  & \mapsto & \sum \beta_i \alpha_i
\end{eqnarray*}
o\`u le complexe   $ \mathcal{P}^{ \cdot} (Y,Z) \otimes \mathcal{P}^{ \cdot} (X,Y)    $ est d\'efini par:
$$
(\mathcal{P}^{ \cdot} (Y,Z) \otimes  \mathcal{P}^{ \cdot} (X,Y))^i = \bigoplus_j \mathcal{P}^{ j} (Y,Z) \otimes  \mathcal{P}^{ i-j} (X,Y), \quad \forall i \in \mathbb Z
$$
avec les diff\'erentielles
$$
  d( \beta   \otimes  \alpha ):= d(\beta) \otimes \alpha + ( -1)^{ \mid \beta \mid} \beta \otimes d(\alpha)
$$

\begin{deft}
\label{dgdp}
Si $\mathcal{P}^{\cdot}$ une $dg-$cat\'egorie, on d\'efinit $ \mathcal{P}^{(\infty,1)}$ une cat\'egorie simplicial telle que:
\begin{itemize}
\item $ Ob(\mathcal{P}^{ \cdot }) = Ob(\mathcal{P}^{( \infty, 1)})  $
\item $ \forall X, Y \in Ob(\mathcal{P}^{(\infty, 1)}) :$
\begin{equation}
\label{pdp}
 \mathcal{P}^{(\infty,1)}(X,Y):= DP( \tau_{\leq 0} \mathcal{P}^{ \cdot}(X,Y) )
\end{equation}
\end{itemize}
\end{deft}

o\`u  $DP( \tau_{\leq 0} \mathcal{P}^{ \cdot}(X,Y) )$ est l'ensemble simplicial de Dold-Puppe \cite{SIMP2} associ\'e au complexe  $\tau_{\leq 0} \mathcal{P}^{ \cdot}(X,Y)$ qui est
 la troncation n\'egative du complexe $\mathcal{P}^{ \cdot}(X,Y)$ d\'efinie par

\[
 (\tau_{\leq 0} \mathcal{P}(X,Y) )^i = \left\{
          \begin{array}{ll}
          \mathcal{P}(X,Y) ^i & \qquad \mathrm{si} \quad i < 0 \\
            Z^0 = Ker( d : \mathcal{P}(X,Y) ^0 \to \mathcal{P}(X,Y) ^1 ) & \qquad \mathrm{si} \quad  i=0\\
               0 & \qquad \mathrm{si}\quad i > 0
          \end{array}
        \right.
\]
et on peut l'\'ecrire:
$$
\tau_{\leq 0} \mathcal{P}^{ \cdot}(X,Y)  = ( \quad ... \rightarrow   \mathcal{P}(X,Y) ^{-i} \rightarrow  ... \rightarrow  \mathcal{P}(X,Y) ^{-2} \rightarrow 
 \mathcal{P}(X,Y) ^{-1}  \rightarrow  Z^0 \quad )
$$

Les fl\`eches de $ \mathcal{P}^{ \cdot }$ se composent avec un produit associatif 

\begin{equation}
\label{assop}
 \forall X,Y,Z \in Ob(\mathcal{P}) : \quad  \mathcal{P}^{(\infty,1)}(Y,Z) \times  \mathcal{P}^{(\infty,1)}(X,Y)   \rightarrow \mathcal{P}^{(\infty,1)}(X,Z)
\end{equation}

Dans le papier \cite{HIRS} on a un morphisme
\begin{equation}
\label{dphisi}
DP(\tau_{ \leq 0}\mathcal{P}^{\cdot}(Y,Z)) \times  DP(\tau_{ \leq 0}\mathcal{P}^{\cdot}(X,Y))   \rightarrow  DP(\tau_{ \leq 0}\mathcal{P}^{\cdot}(Y,Z) ) \otimes 
 DP(\tau_{ \leq 0}\mathcal{P}^{\cdot}(X,Y))
\end{equation}
D'autre part on a un morphisme 
$$
\tau_{\leq 0} \mathcal{P}^{\cdot}(Y,Z) \otimes \tau_{\leq 0} \mathcal{P}^{\cdot}(X,Y) \rightarrow 
\tau_{\leq 0} \mathcal{P}^{\cdot}(X,Z) \otimes \mathcal{P}^{\cdot}(X,Y)  
$$
et on composant avec (\ref{dphisi}) on obtient une multiplication comme dans (\ref{assop})

\begin{equation}
\label{dppasso}
  \forall X,Y,Z \in Ob(\mathcal{P}) : \quad DP(\tau_{ \leq 0}\mathcal{P}^{\cdot}(Y,Z))  \times  DP(\tau_{ \leq 0}\mathcal{P}^{\cdot}(X,Y))   \rightarrow DP(\tau_{ \leq 0}\mathcal{P}^{\cdot}(X,Z))
\end{equation}
Ceci est la composition dans $\mathcal{P}^{(\infty,1)}$.
%%%%%%%%%%%%%%%%%%%%%%%%%%%%%%%%%%%%%%%%%%%%%%%%
% PAGE 2
%%%%%%%%%%%%%%%%%%%%%%%%%%%%%%%%%%%%%%%%%%%%%%%%
\begin{deft}
 Soit $\mathcal{A} = \mathcal{A}^{ (\infty, 1)} $ une $ ( \infty, 1) $-cat\'egorie, on d\'efinit $\mathcal{A}^{( \infty, 0)} $ 
comme la $(\infty,0)$-cat\'egorie int\'erieur \`a $\mathcal{A}$.

Si on consid\`ere $\mathcal{A}^{( \infty, 0)}(X,Y)  $ comme le sous-ensemble simplicial de $\mathcal{A}^{( \infty, 1)}(X,Y)$ en ne prenant que les morphismes inversibles, la sous-cat\'egorie $\mathcal{A}^{( \infty, 0)} \subseteq \mathcal{A}^{( \infty, 1)} $ peut \^etre vue
 comme une quasi-cat\'egorie ou un foncteur vers un ensemble simplicial de Kan, et on note
$$
\mathcal A^{es} := \mathcal{A}^{( \infty, 0)}
$$

\end{deft}

On d\'efinit plus pr\'ecisement le nerf bisimpliciale de cette cat\'egorie par

\begin{eqnarray*}
 (\mathcal{N} \mathcal{A})_{0,m} & = & ob ( \mathcal{A}) \\
(\mathcal{N} \mathcal{A})_{n,m} & = &  \coprod_{\substack {x_0, ..., x_n  \\
x_i \in ob({\mathcal A}) }} 
{\mathcal A}(x_0,x_1)_m  \times ... \times {\mathcal A}(x_{n-1},x_n)_m.
\end{eqnarray*}

Un tel objet simplicial correspond \`a un morphisme

\begin{eqnarray*}
 \Delta^{ \cdot} \times \Delta^{ \cdot} & \to & \mathcal Ens\nonumber \\
(n,m) & \mapsto & (\mathcal{N} \mathcal{A})_{n,m}  \nonumber
\end{eqnarray*}

tel que les $ (\mathcal{N} \mathcal{A})_{n,\cdot} \in \mathcal{E}ns^{\Delta^{\cdot}}$

Les morphismes de S\'egal sont des isomorphismes 

$$
(\mathcal{N} \mathcal{A})_{n,\cdot}  \xrightarrow{\simeq} (\mathcal{N} \mathcal{A})_{ 1,\cdot}
 \times_{(\mathcal{N} \mathcal{A})_{ 0,\cdot} }...\times_{(\mathcal{N} \mathcal{A})_{ 0,\cdot} } (\mathcal{N} \mathcal{A})_{ 1,\cdot}
$$

Pour une $(\infty,1)$-cat\'egorie, on a son int\'erieur
$$
\mathcal{A}^{(\infty,0)} \subset  \mathcal{A}^{(\infty,1)} \quad \text{et} \qquad (\mathcal A^{es})_n := (\mathcal{A}^{( \infty, 0)})_{(n,n)}
$$
ce sont les \'el\'ements de la diagonale.
%%%%%%%%%%%%%%%%%%%%%%%%%%%%%%%%%%%%%%%%%%%%%%%%%%%%%%%%%%%%%%%%%%%%%%%%%%%%

%%%%%%%%%%%%%%%%%%%%%%%%%%%%%%%%%%%%%%%%%%%%%%%%%%%%%%%%%%%%%%%%%%%%%%%%%%%%%%%%%%%%%%%%%%%%%%%%%%%%%%%%%%%%%%%%%%%%%%%%%%%%%%%%%%%%%%%%%%%%%%%

\section*{ \'El\'ements de Maurer-Cartan}

\begin{deft}
 Un \'{e}l\'{e}ment de Maurer-Cartan (M-C) pour $E \in Ob(\mathcal{P})$ est un \'{e}l\'{e}ment $\eta \in \mathcal{P}^1(E,E)$ v\'{e}rifiant l'\'{e}quation de {\em Courbure } suivante
\begin{equation}
 \label{M-C}
\delta (\eta) + \eta . \eta = 0
\end{equation}
on d\'efinit l'ensemble $Ob(MC({\mathcal P}))$ comme l'ensemble des couples $( E, \eta )$ o\`u $E \in Ob(\mathcal{P})$ et $\eta$ est un \'el\'ement de Maurer-Cartan. Un tel couple sera appel\'e $MC$-objet
\end{deft}
Si $(E,\eta) $ et $(F, \zeta)$ sont deux $MC$-objets, on d\'efinit la diff\'erentielle
$$
d_{\eta, \zeta}: \mathcal{P}^i(E,F) \rightarrow \mathcal{P}^{i+1}(E,F)
$$
Par $d_{\eta, \zeta}(a):= da + \zeta .a - (-1)^{i} a.\eta $

\begin{lemme}
\label{lemmenot}
On a $d_{\eta, \zeta}^2=0$ et si on pose
$$
MC({\mathcal P})((E,\eta);(F, \zeta)) := (\mathcal{P}^{\cdot}(E,F), d_{\eta, \zeta})
$$
avec la m\^eme multiplication et identit\'e que $\mathcal{P}$, on obtient une dg-cat\'egorie $MC({\mathcal P})$.
\end{lemme}
%%%%%%%%%%%%%%%%%%%%%%%%%%%%%%%%%%%%%%%%%%%%%%%%%%%%%%%%%%%%%%
% PAGE 3
%%%%%%%%%%%%%%%%%%%%%%%%%%%%%%%%%%%%%%%%%%%%%%%%%%%%%%%%%%%
%lemme 2 !!!
\begin{deft}
 On d\'efinit la cat\'egorie de Maurer-Cartan de ${\mathcal P}$ par la dg-cat\'egorie $MC({\mathcal P})$ avec
\begin{itemize}
 \item $ ob(MC({\mathcal P})) = \{  ( X, \rho); \quad  X \in ob(\mathcal{P}), \quad \rho \in \mathcal{P}^1(X,X)$ avec $  d(\rho) + \rho^2 = 0   \} $
\item $ MC({\mathcal P})((X,\rho),(Y,\mu)) := (\mathcal{P} (X,Y); d_{\rho \mu} )$ et $ d_{\rho \mu}(\cdot) = d(*) +  * \circ \rho + \mu \circ * $. 
\end{itemize}

\end{deft}
%%%%%%%%%%%%%%%%%%%%%%%%%%%%%%%%%%%%%%%%%%%%%%%%%%%%%%%%%%%%%%%%%%%%%%%%%%%%%%%%%%%%%%%%%%%%%%%%%%%%%%%%%%%%%%%%%%%%%%%%%%%%%%%%%%%%%%%%%%%%%%%%
\section*{Hypoth\`eses}
Soit $k$ un corps alg\'ebriquement clos de caract\'eristique nul.
On fixera par la suite une $dg$-cat\'egorie $k$-lin\'eaire $\mathcal P$ qui satisfait aux hypoth\`eses suivant::
\begin{enumerate}
 \item L'ensemble des objets $Ob(\mathcal P)$ est fini.
\item Pour tout $E,F \in Ob(\mathcal P)$, pour tout $i \in \mathbb Z, \quad \mathcal P^i(E,F) $ est un $k$-espace vectoriel de dimension fini.
\item Il existe un indice $n > 0$ tel que pour tout $i < -n$, le $k$-espace vectoriel  $\mathcal P^i(E,F)=0$. 
\end{enumerate}
\begin{rmq}
 Si $\mathcal P^{\cdot}(E,F)$ est un complexe strictement parfait sur $k$, alors les hypoth\`eses $2$ et $3$ sont assur\'ees, et pour l'hypoth\`ese $3$ on peut 
juste dire qu'il existe un $n$ tel que $\mathcal P^i(E,F)=0 $, pour $\mid i \mid > n $ (c'est-\`a-dire: $i \notin [-n, n]$). 
\end{rmq}

\section*{ Le $\mathcal{MC}$-pr\'echamps }
Soit $k$ un corps.
On fixe la dg-cat\'egorie $k$-lin\'eaire $\mathcal{P}$ qui satisfait aux hypoth\`eses ci-dessus. 

Pour tout $k$-alg\`ebre $R$, on d\'efinit une $R-dg$-cat\'egorie $\mathcal{P} \otimes_k R$ telle que 
$$
Ob(\mathcal{P} \otimes_k R) = ob(\mathcal{P})
$$
 et
 $$
   \forall E, F \in ob(\mathcal{P} \otimes_k R) : \quad (\mathcal{P} \otimes_k R)(E,F) : = \mathcal{P}(E,F) \otimes_k R
$$ 
On pose
$$
\mathcal{MC_P}^{dg}(R) := MC(\mathcal{P} \otimes_k R)
$$
$$
\mathcal{MC_P}^{( \infty,1)}(R) : = [ \mathcal{MC_P}^{dg}(R)]^{(\infty,1)}
$$
$$
\mathcal{MC_P}^{es}(R) : = [\mathcal{MC_P}^{(\infty,1)}(R)]^{es}
$$
et on d\'efinit $ \mathcal{MC _P} $ comme le $\infty$-champs associ\'e \`a l'$\infty$-pr\'echamps $\mathcal{MC _P}^{es} $.

On note que la condition $3$ de l'hypoth\`ese, qui s'applique aussi aux $\mathcal{P} \otimes_k R$, implique que $\mathcal{MC_P}^{( \infty,1)}$ est en fait un $(n+1,1)$-pr\'echamps, 
donc $\mathcal{MC_P}$ est un $(n+1)$-champs de $(n+1)$-groupoides.   

%%%%%%%%%%%%%%%%%%%%%%%%%%%%%%%%%%%%%%%%%%%%%%%%%%%%%%%%%%%%%%%%%%%%%%%%%%%%%%%%%%%%%%%%%%%%%%%%%%%%%%%%%%%%%%%%%%%%%%%%%%%%%%%%%%%%%%%%%%%%%%%%

\section*{Lissit\'e formelle}
%%%%%%%%%%%%%%%%%%%
Soit $B$ un anneau commutatif et $I$ un id\'eal de $B$ tel que $I^2=0$. 
Supposons que $(F, \zeta_I)$ est un objet de la dg-cat\'egorie $MC(\mathcal{P} \otimes_k B/I)$ et $(E, \eta)$ un objet de la dg-cat\'egorie 
$MC(\mathcal{P} \otimes_k B)$.

Consid\'erons deux morphismes de MC-objets  $a_I$ et $ \alpha_I$ comme suivant:
$$
(F, \zeta_I) \quad :  \quad \overset{  \alpha_I  }{\underset{  a_I }{ \rightleftarrows}} \quad : \quad  ( E,\eta_I := \eta \otimes _B B/I )
$$
tels que $ a_I$ et $ \alpha_I$ sont inverses dans$H^0(MC(\mathcal{P} \otimes_k B/I)) $.
%%%%%%%%%%%%%%%%%%%%%%%%%%%%%%%%%%%%%%%%%%%%%%%%%%%%%%%%%%%%%%%%%%%%

%%%%%%%%%%%%%%%%%%%%%%%%%%%%%%%%%%%%%%%%%%%%%%%%%%%%
On choisit un rel\`{e}vement (modulo $I$) :
$$
(F, \zeta) \quad :  \quad \overset{  \alpha  }{\underset{  a }{ \rightleftarrows}} \quad : \quad  ( E, \eta)
$$
tels que   $\zeta_I := \zeta \otimes _B B/I   $
 et o\`u
 $a$ et $\alpha$  se r\'eduisent \`a   $a_I$ et $ \alpha_I$ modulo $I$ et les compositions $ \alpha . a$ et $ a. \alpha $ sont donn\'{e}es par:
\begin{equation}
\label{aalpha}
\alpha . a = 1 + d_{\eta, \eta}(g) + u = 1 + \eta g +g \eta +u ;   \quad u \in \mathcal{P}^0(E,E).I
\end{equation}
\begin{equation}
\label{alphaa}
 a.\alpha = 1 + d_{\zeta, \zeta}(h) + v = 1 + \zeta h + h \zeta + v ;  \quad v \in \mathcal{P}^0(F,F).I
\end{equation}
On note que $(F, \zeta)$ sera  un objet de la dg-cat\'egorie $MC(\mathcal{P} \otimes B)$,
si
$$
d(\zeta) + \zeta^2=0,
$$
mais nous ne pouvons pas savoir cela {\em a priori}. On cherche donc \`a modifier $\zeta$.

On pose
\begin{equation}
\label{theta}
 \theta =\zeta + \varepsilon ;  \qquad \varepsilon^2 = 0
\end{equation}
\begin{equation}
\label{zeta}
 d(\zeta) + \zeta^2 = \varphi \in P^2(F,F).I %\qquad ( \zeta^2 \neq 0 \quad et \quad \varphi^2 = 0)
\end{equation}
d'apr\'es \ref{theta} on trouve
\begin{eqnarray*}
 \label{dtheta}
d(\theta) + \theta^2 &  = &  d(\zeta) + \zeta^2+ d(\varepsilon) + \zeta \varepsilon + \varepsilon \zeta \\
 & = & \varphi + d(\varepsilon) + \zeta \varepsilon + \varepsilon \zeta \\
 &=&  \varphi + d_{ \zeta \zeta} ( \varepsilon)
\end{eqnarray*}

On a $ d_{ \zeta_I \eta_I} ( \alpha_I) = 0$, on pose $ \gamma =: d(\alpha) + \eta \alpha - \alpha \zeta \in \mathcal{P}^1(F,E).I $,
on applique $d_{ \zeta_I \eta_I}$ sur $ \gamma$
\begin{eqnarray*}
 d_{ \zeta_I \eta_I} ( \gamma )& = & d( \gamma) + \eta \gamma + \gamma \zeta \\
    &  = &   d( d(\alpha) + \eta \alpha - \alpha \zeta) + \eta (d(\alpha) + \eta \alpha - \alpha \zeta) + (d(\alpha) + \eta \alpha - \alpha \zeta) \zeta \\
   & = &   d^2( \alpha) + d(\eta \alpha) - d( \alpha \zeta) + \eta d(\alpha) + \eta^2 \alpha - \eta \alpha \zeta + d(\alpha) \zeta + \eta \alpha \zeta - \alpha \zeta^2 \\
   & = & d(\eta) \alpha - \eta d(\alpha) - d(\alpha)\zeta - \alpha d(\zeta) + \eta d(\alpha) + \eta^2 \alpha - \eta \alpha \zeta + d(\alpha) \zeta + \eta \alpha \zeta - \alpha \zeta^2   \\
   & = &  ( d( \eta) + \eta^2) \alpha - \alpha ( d(\zeta ) + \zeta^2) \\
   & = & - \alpha \varphi
\end{eqnarray*}
%%%%%%%%%%%%%%%%%%%%%%%%%%%%%%%%%%%%%%%%%%%%%%%%%%%%%%%%%%%%%%%%%%%%%%%%
% PAGE 4
%%%%%%%%%%%%%%%%%%%%%%%%%%%%%%%%%%%%%%%%%%%%%%%%%%%%%%%%%%%%%%%%%%%%%%%%
Donc on a
\begin{equation}
 \label{dgamma}
d_{ \zeta_I \eta_I}(\gamma) = - \alpha_I \varphi
\end{equation}
 on calcule aussi
\begin{eqnarray*}
 d_{ \zeta_I \zeta_I} (a_I \gamma )& =  &   d_{ \eta_I \zeta_I} (a_I) \gamma + a_I  d_{ \zeta_I \zeta_I} ( \gamma)\\
  & = &   -a_I \alpha_I \varphi  \\
  & = & - ( 1 +  d_{ \zeta \zeta}(h)) \varphi
\end{eqnarray*}
or l'identit\'e de Bianchi montre que $ d_{\zeta \zeta}( \varphi) = 0$
et on a
\begin{eqnarray*}
 d_{ \zeta \zeta} (h \varphi )& =  &   d_{ \zeta \zeta}(h) \varphi + h d_{ \zeta \zeta}( \varphi)\\
  & = &   d_{ \zeta \zeta}(h) \varphi
  \end{eqnarray*}
donc
\begin{eqnarray*}
 d_{ \zeta_I \zeta_I} (a_I \gamma )& =  &   d_{ \eta_I \zeta_I} (a_I) \gamma + a_I  d_{ \zeta_I \zeta_I} ( \gamma)\\
   & = & - ( 1 +  d_{ \zeta \zeta}(h)) \varphi  \\
   & = & - \varphi - d_{ \zeta \zeta}(h) \varphi \\
   & = & - \varphi - d_{ \zeta \zeta} (h \varphi )
\end{eqnarray*}
et donc
$$
\varphi = - ( d_{ \zeta \zeta}( a \gamma)+ d_{ \zeta \zeta} (h \varphi) )
$$
dans ce cas, on choisit $ \varepsilon = a \gamma + h \varphi $.

Etape2 :
On rempla\c cant $\zeta$ par $ \theta$, nous pouvons supposer que $d(\zeta) + \zeta^2 =0$, par cons\'equent $d_{\zeta \eta}^2(\varepsilon)=0$
%%%%%%%%%%%%%%%%%%%%%%%%%%%%

On modifie $\zeta$ et $ \alpha$ de sorte que $ d_{ \zeta \eta} ( \alpha) = 0 $ et $ d_{\zeta \eta }^2=0 $.

On note par $ \omega : = d_{  \zeta \eta } ( \alpha) = d(\alpha) + \eta \alpha - \alpha \zeta $ alors $ d_{\zeta \eta} ( \omega) = 0$

On pose $ \theta = \zeta + \varepsilon'$ pour un $ \varepsilon' $ choisi de sorte que $ d_{ \zeta \zeta}(\varepsilon') = 0$, donc on a bien
$$
 d(\theta) + \theta^2 = d(\zeta) + \zeta^2 + d_{ \zeta \zeta} ( \varepsilon') = 0
$$
donc
\begin{eqnarray*}
 d_{  \theta \eta }(\alpha)  & =  & d_{ ( \zeta + \varepsilon' )\eta } (\alpha) \\
                             & = &  d( \alpha)+ \eta \alpha +  \alpha ( \zeta  + \varepsilon' )  \\
                             & = & d( \alpha) +  \eta \alpha + \alpha \zeta + \alpha \varepsilon' \\
                             & = & (d(\alpha) + \eta \alpha + \alpha \zeta ) + \alpha \varepsilon' \\
                             & = & d_{\zeta \eta}(\alpha) + \alpha \varepsilon' 
 \end{eqnarray*}

Prenant maintenant $ t \in \mathcal{P}^1 (F,E) .I $, et on calcule
\begin{eqnarray*}
 d_{\theta \eta } ( \alpha + t)  & = &   d_{\theta \eta}  ( \alpha)  +  d_{\theta \eta}  (t)  \\
                                 & = & d(\alpha) + \eta \alpha + \alpha \theta + d(t) + \eta t + t \theta   \\
                                 & = &  d(\alpha) + \eta \alpha + \alpha ( \zeta + \varepsilon')  + d(t) + \eta  t + t  ( \zeta + \varepsilon')  \\
                                 & = & d(\alpha) + \eta \alpha + \alpha \zeta + \alpha \varepsilon'   + d(t) + \eta t +  t \zeta + t \varepsilon'   \\
                                 & = & (d(\alpha) + \eta \alpha + \alpha \zeta ) +\alpha  \varepsilon'  + (d(t) + \eta t  + t \zeta ) + t \varepsilon' \\
                                 & = & d_{\zeta \eta}(\alpha) + d_{\zeta \eta}(t) +\alpha \varepsilon'  \qquad (t \varepsilon'  = 0 mod (I)) \\
\end{eqnarray*}
Posons maintenant $\varepsilon'= - a \omega$ donc
\begin{eqnarray*}
 \alpha \varepsilon' & = & - \alpha a \omega \\
                     & = & - (1+ d_{\eta \eta}(g)) \omega \\
                     & = & - \omega - d_{\eta \eta}(g) \omega \\
                     & = & - \omega - d_{\zeta \eta}(g \omega)
\end{eqnarray*}
car
\begin{eqnarray*}
 d_{\zeta \eta}(g \omega) & = & d_{\eta \eta}(g)\omega + g d_{\zeta \eta}(\omega)\\
                         & = &    d_{\eta \eta}(g)\omega \qquad ( car \quad  d_{\zeta \eta}(\omega) =0).
\end{eqnarray*}

Posons $t= g \omega$ et rempla\c cons \c ca dans la formule de $d_{\theta \eta} (\alpha +t)$, on obtient
\begin{eqnarray*}
 d_{\theta \eta} (\alpha +t)   & = & d_{\zeta \eta }(\alpha) + d_{\zeta \eta }( g \omega)  + \alpha \varepsilon' \\
                               & = & d_{\zeta \eta}(\alpha) +  d_{\zeta \eta }( g \omega)  - \omega - d_{\zeta \eta }( g \omega)\\
                               & = & d_{\zeta \eta}(\alpha) - \omega \\
                               & = & 0            
\end{eqnarray*}
 car par hypoth\`ese $ d_{\zeta \eta}(\alpha) = \omega$.

Ces r\'esultats montrent bien le th\'eor\`eme suivant:
\begin{thm}
\label{thmrelevement}
 Il existe des relevements $\zeta$ pour $\zeta_I$ et $\alpha$ pour $\alpha_I$ tels que $d(\zeta) + \zeta^2=0$ et $d_{\zeta \eta}(\alpha)=0$.
\end{thm}

%%%%%%%%%%%%%%%%%%%%%%%%%%%%%%%%%%%%%%%%%%%%%%%%%%%%%%%

\section*{Construction de la carte $ \varphi : V \to \mathcal{MC_P}$}

On choisit des $k$-bases $  e^1, \ldots, e^r   \in \mathcal{P}^{1}(E,E) $ et $  f^1, \ldots, f^s    \in \mathcal{P}^2(E,E) $, et on d\'efinit l'application 
\begin{eqnarray}
\label{prodp1}
 \mathcal{P}^{1}(E,E) \times \mathcal{P}^{1}(E,E) & \to & \mathcal{P}^{2}(E,E) \\
(e^i,e^j) & \mapsto & e^ie^j = \sum_{l=1}^s a^{ij}_lf^l \nonumber
\end{eqnarray}
avec
\begin{equation}
 \label{p1}
\mathcal{P}^{1}(E,E) = k^r = \{ \sum x_ie^i, \quad x_i \in k, e^i \in \mathcal{P}^{1}   \}
\end{equation}
et
$$
\mathcal{P}^{2}(E,E) = k^s = \{ \sum y_lf^l, \quad y_l \in k, f^l \in \mathcal{P}^{2}   \}
$$
les coefficients structurels  $a^{ij}_l \in k$ sont unique. 
La diff\'erentielle 
 \begin{eqnarray}
 d :\mathcal{P}^{1}(E,E)  & \to & \mathcal{P}^{2}(E,E) \\
e^i & \mapsto & d(e^i) = \sum_{l=1}^s y^i_lf^l \nonumber
\end{eqnarray}
et $ (z^1, \ldots, z^s) \in k^s $.
%%%%%%%%%%%%%%%%%%%%%%%%%%%%%%%%%%%%%%%%%%%%%%%%%%%%%%%
Pour un $\eta = \sum_{i=1}^rz_ie^i  \in \mathcal{P}^{1}(E,E)$, on a
\begin{eqnarray*}
 \eta^2  & = & (\sum_{i=1}^rz_ie^i)(\sum_{j=1}^rz_je^j)  \\
         & = &  \sum_{i,j}z_iz_je^ie^j                     \\
         & = & \sum_{i,j}z_iz_j ( \sum_{l=1}^s a^{ij}_lf^l) \\
         & = &  \sum_{l=1}^s (\sum_{i,j}z_iz_j  a^{ij}_l)f^l.
\end{eqnarray*}
et 
\begin{eqnarray*}
 d(\eta) & = & d(  \sum_{i=1}^rz_ie^i) \\
         & = & \sum_{i=1}^r z_i d(e^i) \\
         & = & \sum_{i=1}^r z_i (\sum_{l=1}^s b^i_lf^l) \\
         & = & \sum_{l=1}^s (\sum_{i=1}^r z_i  b^i_l)f^l
\end{eqnarray*}
Donc
\begin{eqnarray*}
 d(\eta) + \eta^2 &  =  &  \sum_{i=1}^r z_i (\sum_{l=1}^s y^i_lf^l) + \sum_{i,j}z_iz_j ( \sum_{l=1}^s a^{ij}_lf^l) \\
                  &  =  &  \sum_{l=1}^s ( \sum_{i=1}^r z_i y^i_l)f^l +  \sum_{l=1}^s(\sum_{i,j}z_iz_j  a^{ij}_l)f^l \\
                  &  =  &  \sum_{l=1}^s  (  \sum_{i=1}^r z_i y^i_l +\sum_{i,j}z_iz_j  a^{ij}_l)f^l \\
                  &  =  &  \sum_{l=1}^s  c_l  f^l \quad et \quad c_l = \sum_{i=1}^r z_ib^i_l +\sum_{i,j}z_iz_j  a^{ij}_l
\end{eqnarray*}
Notre question est de construire un foncteur $ Alg-Com_k \to \mathcal{E}ns$ qui a pour chaque anneau ( ou alg\'ebre) $B$ son image l'nsemble des points de $B$  et montrer qu'il est repr\'esentable par un sch\'ema affine qu'on notera par $V_E$. Pour cela 
on d\'efinit les sch\'emas affines par le foncteur $Spec$ et on les note par 
$$
\mathcal P^1_{sch}(E,E): = Spec( k[x_1, \ldots, x_r]  )
$$
et
$$
\mathcal P^2_{sch}(E,E): = Spec( k[y_1, \ldots, y_s]  )
$$

%%%%la place de $\mathcal{P}^{1}(E,E)  $ d\'efinit par le formule $(\ref{p1})$ on choisit le sch\'ema affine
$$
\mathcal{P}^{1}_{sch}(B) = \mathcal{P}^{1}(E,E) \otimes_kB = \mathbb{A}^r(B) = B^r
$$
On peut ensuite utiliser les m\^emes \'equations pour des $z_i \in B$ et $B$ un alg\`ebre commutatif et on d\'efinit l'espace $B^r$ par
\begin{eqnarray*}
 B^r  & = & k^r \otimes_kB \\
      & = & \{ (z_1, \ldots, z_r) : z_i \in B \quad \forall i \in \{1, \ldots, r  \}   \} \\
      & = & \{ \sum_{i=1}^r e^i \otimes z_i, \quad z_i \in B, e^i \in  \mathcal{P}^{1}(E,E) \} \\
      & = & \mathcal{P}^{1}(E,E) \otimes_kB
\end{eqnarray*}
Comme
$$
(\mathcal{P}^{1}(E,E) \otimes_kB) \otimes_k (\mathcal{P}^{1}(E,E) \otimes_kB) =  (\mathcal{P}^{1}(E,E)  \otimes_k \mathcal{P}^{1}(E,E)) \otimes_kB
$$
On peut d\'efinir un produit similaire que $(\ref{prodp1})$ par:

\begin{eqnarray*}
 (\mathcal{P}^{1}(E,E)  \otimes_k \mathcal{P}^{1}(E,E)) \otimes_kB & \to &  \mathcal{P}^{2}(E,E) \otimes_kB \\
(\sum_{i=1}^r e^i \otimes z_i, \sum_{j=1}^r e^j \otimes z_j) & \mapsto & \sum_{i,j=1}^r e^ie^j \otimes z_iz_j
\end{eqnarray*}
d'apr\`es la formule $(\ref{prodp1})$  pour tout $ \eta_B \in \mathcal{P}^{1}(E,E) \otimes_kB$ on peut  calculer $ \eta^2 \in  \mathcal{P}^{2}(E,E) \otimes_kB$ par 
\begin{eqnarray*}
\eta_B^2   & = &   (\sum_{i=1}^r e^i \otimes z_i )( \sum_{j=1}^r e^j \otimes z_j) \\
                    & = &  \sum_{i,j=1}^r e^ie^j \otimes z_iz_j \\
                    & = &   \sum_{i,j=1}^r (\sum_{l=1}^s a^{ij}_lf^l) \otimes z_iz_j  \\  
                    & = &  \sum_{l=1}^s (  \sum_{i,j=1}^r  a^{ij}_l\otimes z_iz_j ) f^l  \\                             
\end{eqnarray*}
et aussi
\begin{eqnarray*}
d(\eta_B)           & = &   d(\sum_{i=1}^r e^i \otimes z_i ) \\
                    & = &  \sum_{i=1}^r d(e^i \otimes z_i)  \\
                    & = &      \sum_{i=1}^r (d( e^i) \otimes z_i + e^i \otimes d(z_i)) \\
                    & = &  \sum_{i=1}^r d( e^i) \otimes z_i  \qquad car \quad d(z_i) = 0 \quad \forall i = 1, \ldots, r \\
                    & = &     \sum_{i=1}^r ( \sum_{l=1}^s b^i_lf^l) \otimes z_i    \\
                    & = &    \sum_{l=1}^s ( \sum_{i=1}^r  b^i_l \otimes z_i ) f^l   \\              
\end{eqnarray*}
Donc
$$
 d(\eta_B) + \eta_B^2 =  \sum_{l=1}^s ( \sum_{i=1}^r  b^i_l \otimes z_i ) f^l  + \sum_{l=1}^s (  \sum_{i,j=1}^r  a^{ij}_l\otimes z_iz_j ) f^l\\
$$

%%%%%%%%%%%%%%%%%%%%%%%%%%%%%%%%%%%%%%%%%%
%%%%%%%%%%%%%%%%%%%%%%%%%%%%%%%%%%%%%%%%%%         je suis la!!
On note le foncteur $V_E$ par
\begin{equation}
 \label{cartevmc}
V_E(B) = \{ \sigma \in \mathcal{P}^1(E,E) \otimes_k B, \quad \mbox{tel que } \quad d(\sigma) + \sigma^2 = 0  \}
\end{equation}
l'ensemble des \'el\'ements de Maurer-Cartan de $ \mathcal{P}^1(E,E) \otimes_k B $ qui peut  \^etre vue comme un foncteur repr\'esentable par un sch\'ema affine
$$
V_E(B) : k-Alg \to \mathcal{E}ns
$$
%%%%%%%%%%%%%%%%%%%%%%%%%%%%%%%%%%%%%%%%%%%%%%%%%%%%%%%%%%%%%%%%%%%%
% PAGE 5
%%%%%%%%%%%%%%%%%%%%%%%%%%%%%%%%%%%%%%%%%%%%%%%%%%%%%%%%%%%%%%%%%%%

\begin{proposition}
 Le foncteur $V_E(B)$ est le pr\'eimage de l'origine de l'epace affine $\mathcal{P}^2_{sch}(E,E)$, par le morphisme $ {\rm courb} $.
 Par cons\'equent, $V_E(B)$ est repr\'esent\'e par le sch\'ema affine
$$
V_E(B) = Spec \left( \frac{k[x_1, \ldots, x_r]}{\mathcal{I}} \right)
$$
o\`u $ \mathcal{I} $ est l'id\'eal engendr\'e par les \'el\'ements 
$$
{\rm courb}^{\ast}(y_l) =  \sum_{i=1}^r b^i_lx_i + \sum_{i=1}^r \sum_{j=1}^r a^{ij}_l x_ix_j.
$$
\end{proposition}
Nous obtenons donc un sch\'ema affine $ V_E $ et un morphisme $ V_E \to \mathcal{MC_P} $.

\begin{corr}
 Le morphisme $ V_E \to \mathcal{MC_P} $ est formellement lisse. 
\end{corr}
{\bf \em Preuve:}
Voir Th\'eor\`eme \ref{thmrelevement}.

\begin{proposition}
\label{prop2}
 Soit $X,Y$ des sch\'emas avec des morphismes 
$$
f:X \to \mathcal{MC_P}, \qquad g: Y \to \mathcal{MC_P} 
$$
alors le produit fibr\'e 
$$
X \times_{ \mathcal{MC_P}} Y
$$ 
est repr\'esent\'e par un $n$-champ d'Artin (ie : g\'eom\'etrique ) de type fini.
\end{proposition}
{\bf \em Preuve:}
On peut supposer que $X$ et $Y$ sont affines.
Soient $B$ et $C$ deux alg\`ebres.  
On pose $X = Spec(B)$ et $ Y= Spec(C) $ les sch\'emas affines associ\'es, et on peut supposer que les morphismes $f$ et $g$ proviennent 
d'\'el\'ements de $\mathcal{MC_P}^{es}(B)$ et $\mathcal{MC_P}^{es}(C)$ qu'on notera $(F, \eta)$ et $(G,\rho)$.

Ce sont alors des \'el\'ements de Maurer-Cartan $\eta \in \mathcal{P}^1(F,F)\otimes_kB$ et $\rho  \in \mathcal{P}^1(G,G)\otimes_kC$.
Une fl\`eche  $Spec(R) \to X \times Y$ correspond \`a $B \otimes_kC \to R$, et on notera $\eta_R$ et $\rho_R$ les images de $\eta \otimes 1_C$
 et $1_B \otimes \rho$ dans $\mathcal P^1 (F,F)\otimes_kR$ et $\mathcal P^1 (G,G)\otimes_kR$ respectivement.
Ce sont des MC-\'el\'ements.

L'ensemble simplicial des relevements 
$$
\begin{array}{ccc}
   &   &  X \times_{ \mathcal{MC_P}} Y\\
   & \nearrow & \downarrow  \\
Spec(R) & \to & X \times Y 
\end{array}
$$  
est \'equivalent \`a 
$$
\mathcal M^{(\infty,0)}(R) ((F, \eta_R),(G,\rho_R))
$$
qui est \'egale au sous-ensemble simplicial 
$$
\mathcal M^{(\infty,1)}(R) ((F, \eta_R),(G,\rho_R))^{iso} \subset \mathcal M^{(\infty,1)}(R) ((F, \eta_R),(G,\rho_R))
$$
des morphismes inversibles. 

D'autre part 
$$
\mathcal M^{(\infty,1)}(R) ((F, \eta_R),(G,\rho_R)) = DP( \tau_{\leq 0}(\mathcal P(F,G)\otimes_kR, d_{\eta_R, \rho_R} )).
$$
Posons 
$$
L^{\cdot} : = (P(F,G)\otimes_kB\otimes_kC, d_{\eta \otimes 1_C, 1_B \otimes \rho } )
$$
c'est un complexe parfait de $B\otimes_kC$-modules.

Le foncteur 
$$
R \mapsto DP( \tau_{\leq 0}(\mathcal P(F,G)\otimes_kR, d_{\eta_R, \rho_R} )) =  DP( \tau_{\leq 0} ( L^{\cdot} \otimes_{B \otimes_kC} R))
$$
d\'efinit un $n$-champs d'Artin $H(f,g)$ au dessus de $X \times Y$ d'apr\`es \cite{SIMP1}.

Le produit fibr\'e $X \times_{ \mathcal M}Y$ est le sous champs des morphismes inversibles $H(f,g)^{iso} \subset H(f,g)$.
On pourra voir que c'est un  sous-champs ouvert, ce qui montre que  $X \times_{ \mathcal M}Y$ est un $n$-champs d'Artin.

\begin{lemme}
 Soit $k$ un corps alg\'ebriquement clos, $M$ un $(n+1)$-champs sur un sch\'emas affine $k-Alg-Com$ et si $V$ est un sch\'ema de type fini sur $k$, et si on a
un morphisme $\varphi : V \to M.$ 
Si $\varphi$ est formellement lisse, et si l'application $ V(k) \to \mathcal{P}^0(M(k)) $ est surjective, et si en plus de \c ca $M$ satisfait les conditions de 
la proposition (\ref{prop2}). Alors $M$ est un $(n+1)$-champs d'Artin de type fini. 
\end{lemme}

\begin{corr}
\label{mcartin}
 le $\mathcal{MC_P}$ est un $(n+1)-$champs d'Artin de type fini.
\end{corr}
%%%%%%%%%%%%%%%%%%%%%%%%%%%%%%%%

Pridham montre qu'un $(n+1)-$champs d'Artin de type fini peut \^etre repr\'esent\'e par un sch\'ema simplicial ayant des bons propri\'et\'es \cite{PRID}. 
Il serait int\'eressant de construire un tel sch\'ema simplicial explicitement en partant des cartes $V_E$.

\section*{Exemple:}

Dans cet exemple, on va utiliser les m\^emes arguments que dans $\cite{BENZ08}$, pour cela on donne une bref rappel de l'$\infty$-champs $Perf$:

Un complexe strictement parfait est un complexe strictement born\'{e} des fibr\'{e}s vectoriels alg\'{e}briques.
On peut dire qu'un complexe est  parfait s'il est locallement quasi-isomorphe \`a un complexe strictement parfait.  

Le champs $Perf$ des complexes parfait est un foncteur

\begin{eqnarray*}
 Aff_k^o & \to & Ens^{\Delta^o}  \\
     X & \mapsto & \mathcal{N}erf(CpxPerf_{qiso}(X)) 
\end{eqnarray*}
o\`u $ Aff_k^o$ est la cat\'egories des sch\'emas affines et $CpxPerf_{qiso}$ est la cat\'egorie dont les objets sont des complexe parfait et les fl\`eches sont les quasi-isomorphismes entre ces complexes.

 On choisit un nombre fini des complexes d'espaces vectoriels de dimensions finies avec $d=0$, on les notes par $ E^{\cdot} : = (E^{\cdot},0)$.
 Soit $\mathcal P^{\cdot}$ la dg-cat\'egorie form\'ee par ces objets, avec
$$
\mathcal P^{\cdot}(E^{\cdot},F^{\cdot}) : = (E^{\cdot})^* \otimes F^{\cdot}, \quad d=0.
$$
 
 Un \'el\'ement de Maurer-Cartan $\eta$ dans $ \mathcal{P}^1(E^{\cdot},E^{\cdot})$ est un diff\'erentiel 
 $$
 \eta = \{ \eta^i \}, \quad \eta^i : E^i \to E^{i+1} \quad \text{et} \quad \eta^{i+1}\eta^i = 0.
 $$
qui permet d'obtenir un nouveau complexe
$$
\cdots \xrightarrow{\eta^{i-1}} E^i \xrightarrow{\eta^{i}} E^{i+1} \xrightarrow{\eta^{i+1}} \cdots 
$$

Le sch\'ema $V_{E^{\cdot}}$ est la vari\'et\'e des complexes de Buchsbaum-Eisenbud,
voir \cite{BUCH1}, \cite{BUCH2}, \cite{BRUN}, \cite{HUNE}, \cite{KEMP}, \cite{MASS},  \cite{TRIV} et \cite{YOSH}. 
Le champs $\mathcal{MC_P}$ est le sous-$(n+1)$-champs ouvert de $Perf$ couvert par les cartes
$V_{E^{\cdot}}$ et les morphismes 
$$
V_{E^{\cdot}} \to \mathcal{MC_P} \subseteq Perf.
$$ 
sont ceux mentionn\'es dans $\cite{BENZ08}$. 
\\
L'application du corollaire \ref{mcartin} dans le cas de cet exemple, nous donne une nouvelle d\'emonstration du th\'eor\`eme de Toen-Vaquie dans \cite{TOEN2} que le champs $Perf$ est recouvert de $n$-champs d'Artin.

%\newpage

\noindent
{\sc Laboratoire J.A.\ Dieudonn\'e, Universit\'e de Nice-Sophia
Antipolis, Parc Valrose, 06108 Nice Cedex 02, France},
\verb+bbrahim@unice.fr+

\end{document}